\documentclass[11pt,reqno]{amsart}

\usepackage[utf8]{inputenc}
\usepackage[T1]{fontenc}
\usepackage{amsmath,amssymb,amsthm}
\usepackage{graphicx}
\usepackage{hyperref}
\usepackage{lineno}

\usepackage[margin=3cm,footskip=1cm]{geometry}
\usepackage{caption}
\newtheorem*{uconj}{Conjecture 9}
\theoremstyle{definition}

\title[AAALS Helmholtz]{AAA least squares solution of Helmholtz problems}

\author{Stefano Costa}
\thanks{S. Costa: IEEE Member, Piacenza, Italy. Email: \texttt{stefano.costa@ieee.org}}


\date{January 26, 2026}
\begin{document}
	% \linenumbers
	
	\begin{abstract}
		This paper presents an adaptive numerical framework for solving exterior ``sound-soft'' scattering problems governed by the Helmholtz equation. By interpreting the Method of Fundamental Solutions through the lens of rational approximation, we introduce an automated strategy for singularity placement based on the analytic continuation of boundary data. The proposed AAALS-Helmholtz algorithm leverages a ``continuum'' variant of the AAA algorithm to identify the singularities limiting analytic extension, and to ensure an optimal source distribution even for complex, non star-shaped geometries. Furthermore, we establish a formal connection between the Helmholtz and Laplace problems, providing a theoretical justification for the ``double poles'' technique. The approach offers a robust, meshless alternative to heuristic source placement.
	\end{abstract}
	
	\maketitle
	
	\section{Introduction}\label{s:1}
	The numerical solution of boundary value problems (BVPs) for the Laplace and Helmholtz equations is a cornerstone of computational acoustics and scattering theory. Among the various techniques available, the Method of Fundamental Solutions (MFS), also known as the charge simulation method or the method of auxiliary sources, has garnered significant attention due to its conceptual simplicity and meshless nature (\cite{Katsurada88,Rokhlin90,Fairweather98}). The fundamental principle of the MFS is to approximate the solution $u$ by a linear combination of fundamental solutions of the governing equation, with singularities located outside the domain of interest $\Omega$. For the Helmholtz equation,
	\begin{align}
		\nonumber  \Delta u + k^2 u = 0 \quad \text{in } \Omega \\
		u = v \quad \text{on } \partial\Omega
		\label{eq:helmholtz}
	\end{align}
	where $k$ is the wavenumber, the approximation in the literature typically takes the form:
	\begin{equation}
		u(z) \approx u^{(J)}(z) = \sum_{j=1}^{J} \alpha_j H_0^{(1)}(k|z-p_j|),
		\label{eq:helmholtz-approx}
	\end{equation}
	where $H_0^{(1)}$ denotes the Hankel function of the first kind (outward-radiating) and $\{p_j\}$ represents a set of source points located in the exterior region $\mathbb{C} \setminus \overline{\Omega}$. Despite its potential for exponential convergence on analytic domains, the classical MFS suffers from well-documented robustness and accuracy issues: the effective placement of the singularities $\{p_j\}$. The resulting linear systems are often severely ill-conditioned, and the numerical stability of the representation is intrinsically linked to the extent to which the solution can be analytically continued into the complex plane. For domains exhibiting corners or non-convex geometries, manual placement of these singularities becomes increasingly difficult.
		
	Recent advances in rational approximation (\cite{Nakatsukasa18}) have provided a solution to this placement problem for the Laplace equation. The AAA least squares (AAALS) method of \cite{Costa20} and \cite{CostaTrefethen23} interprets the MFS through the lens of rational approximation, where singularity locations are determined adaptively rather than by a priori heuristics, completely changing its efficiency. This approach has demonstrated that rational approximants can automatically resolve singularities near corners, effectively ``upgrading'' the MFS for complex geometries. However, extending this adaptive philosophy to the Helmholtz equation presents distinct challenges. While ``lightning'' solvers have been proposed to cluster poles near corners based on exponential formulas (\cite{GopalTrefethen19,Trefethen20_lightning}), a fully adaptive ``AAALS-Helmholtz'' method for general domains and wave numbers $k\gtrsim20$ has remained an open problem, as discussed in detail in \cite{NakatsukasaTrefethen26}. Unlike the Laplace case, where fundamental solutions are poles or monopoles, the Helmholtz operator requires the manipulation of complex Hankel functions.
		
	This paper formalizes and expands upon the main results of \cite{BarnettBetcke08}, which mainly concentrates on ill-conditioning, as well as a number of preliminary experiments by the author and colleagues, starting from \cite{Trefethen2023Helmholtz}, which suggested that an adaptive AAALS-Helmholtz method is indeed viable. We focus specifically on the ``sound-soft'' scattering problem in the exterior of a domain, calculating the response to an incident plane or point-source wave. In this configuration, singularities are placed in the interior of the scatterer, each corresponding to an outgoing circular wave satisfying the Sommerfeld radiation condition. We address several open questions regarding this methodology, providing the first comprehensive experimental treatment of the AAALS-Helmholtz method for scattering problems with high wave number, and establishing a theoretical connection with AAALS for Laplace problems.

	\section{Analytic Continuation and the Schwarz Function}\label{s:2}

	Throughout this paper, we draw upon several key results from complex analysis concerning analytic continuation. Let $\Gamma$ be a parametric analytic Jordan curve in the complex plane $z \in \mathbb{C}$, and let $v(z)$ be a real-valued analytic function defined on $\Gamma$. 
	According to Millar \cite{Millar80,Millar86}, the analytic continuation of the Helmholtz solution $u(z)$ across such a curve is limited by two distinct sources of singularities. The first originates from the singularities of the analytic continuation of the boundary data $v$ itself. The second is intrinsic to the geometry of the curve, arising specifically from the singularities of the Schwarz function $S(z)$, which satisfies $S(z) = \bar{z}$ on $\Gamma$. 
	Since the Schwarz function effectively tracks the location of singularities across boundaries, the analytic continuation of the solution is fundamentally constrained by the singularities of $S(z)$. The critical role of these geometric singularities in governing the transition from polynomial to accelerated rational convergence rates for Laplace problems has been recently characterized by Trefethen \cite{Trefethen24}. Furthermore, as emphasized by Barnett and Betcke \cite{BarnettBetcke08} in the context of Helmholtz problems, the correct positioning of these singularities outside the domain is crucial to obtaining well-behaved bases for the MFS.
	
	Davis and Pollak investigated the properties of $S(z)$ in \cite{DavisPollak58}, demonstrating its derivation from a composition of mappings. Given any pair $(M,m)$ where $M$ is a generic analytic mapping from the unit circle $C(t) = e^{i\pi t}$ in the $\zeta$-plane to $\Gamma$, and $m$ is its inverse from $\Gamma$ to $C$, the Schwarz function can be expressed as:
	\[S(z) = \overline{M}(1/m(z)). \]
	While the properties of $S$ hold irrespective of the specific choice of $(M,m)$, for our purposes we are particularly interested in the mapping $M(\zeta)$, which maps each point $\zeta \in C(t)$ to the corresponding point $z \in \Gamma(t)$. Given a particular choice for $M$, since both $\Gamma$ and $C$ are analytic, it possesses an analytic continuation beyond $C$. Furthermore, the analyticity of $\Gamma$ implies that $M' \neq 0$ on $C$, ensuring the existence of an annular-like region $A \supset C$, bounded internally and externally by some curves $\gamma$ and $\delta$, where $M' \neq 0$. Consequently, $M$ maps $A$ conformally onto the region $M(A) \supset \Gamma$.
		
\begin{figure}[!t]
		\centering{
			\includegraphics[width=0.9\linewidth]{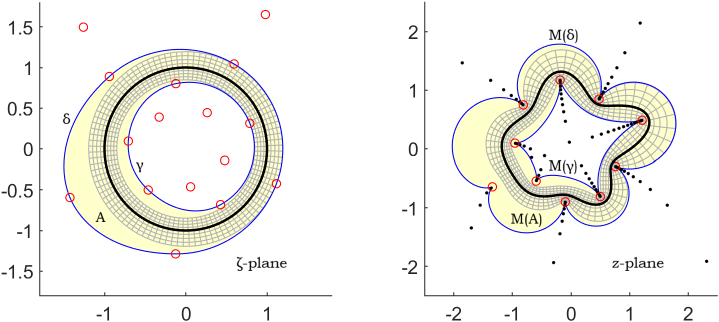}
		}
		\caption{Transformation of the unit circle $C(t)$ in the $\zeta$-plane onto an analytic Jordan curve $\Gamma(t)$ in the $z$-plane (thick black lines) via a generic mapping function $M:\zeta \rightarrow z$. This implicitly defines a conformal map between the regions $A \supset C$ and $M(A) \supset \Gamma$ (shaded yellow). On the left: a circular annular region (gray grid) centered on $C$, the zeros of $M'$ (red circles), and two curves $\gamma$ and $\delta$ (thick blue) connecting the zeros nearest to $C$; these bound a domain of conformality $A$ larger than the annulus, though not necessarily maximal. On the right: their respective images through $M$ and the branch cuts of the Schwarz function for $\Gamma$, delineated by black dots.}
		\label{fig:confmap}
	\end{figure}

	To illustrate this, consider Figure~\ref{fig:confmap}, where a circular annular grid is defined between $C$ and the nearest interior and exterior zeros of $M'$ (the approximant to $M$ has all its poles, not shown, near the origin, as it typically happens). While the curves $\gamma$ and $\delta$ connecting the nearest zeros bound a larger domain of conformality $A$ induced by $M$, determining its maximal analytic extension is computationally challenging; therefore, a circular annular region serves only as a highly conservative approximation. Notably, these nearest zeros are mapped through $M$ to the branch points of the Schwarz function for $\Gamma$. While these geometric observations are established in the literature, the computational methods used to exploit them, discussed in the next section, represent a novel contribution. See \cite{Trefethen25} for a thorough discussion on the numerical computation of the Schwarz function in particular by means of the AAA algorithm.

	\section{Continuum AAA Approximation for Generic Curves}\label{s:3}
	The AAA algorithm begins with a vector $Z$ of real or complex sample points and a vector $F$ of corresponding real or complex function values. Its objective is to find a rational function $r$ such that the $\infty$-norm error $\|r(Z) - F\|$ is minimal; typically, the algorithm aims for a relative error $\|r(Z) - F\| / \|F\| \le 10^{-13}$, achieving close to machine precision in standard 16-digit arithmetic. To ensure numerical stability, the approximation is constructed in barycentric form:
	\begin{equation}
		r(z) = \frac{N(z)}{D(z)} = \frac{\sum_{j=0}^J \dfrac{f_j w_j}{z - s_j}}{\sum_{j=0}^J \dfrac{w_j}{z - s_j}},
	\end{equation}\\
	where $\{s_j\}$ are $J+1$ support points, $\{f_j\}$ are the associated function values, and $\{w_j\}$ are the barycentric weights. The standard AAA algorithm constructs the rational approximant by iteratively populating the set of support points from the larger, predefined set $Z$. A distinctive property of this approach is that the resulting rational function interpolates the data at the support points, ensuring the approximation error vanishes there. In other words, although each $s_j$ appears to be a pole for both the numerator and the denominator, these singularities cancel out in the quotient; thus, each $s_j$ is a removable singularity where $\lim_{z \to s_j} r(z) = f_j$. This barycentric approach is crucial for the success of the algorithm. Traditional representations of rational functions as quotients of polynomials $p(z)/q(z)$ in global bases (like monomials) are notoriously unstable, often requiring extended precision arithmetic to handle functions with poles and zeros clustered near singularities. In contrast, the AAA algorithm remains robust. Furthermore, the barycentric form allows for the efficient extraction of the poles and zeros of the approximant, which can be computed by solving a generalized eigenvalue problem involving the support points and weights. For problems defined on continuous domains, this requires an a priori discretization; however, such a static approach is often inadequate when singularities are located near the boundary. Without manual, high-density clustering of the initial sample set near these regions the available candidates for support points may not suffice to resolve the function's rapid variations.
		
\begin{figure}[!t]
		\centering{
			\includegraphics[width=0.7\linewidth]{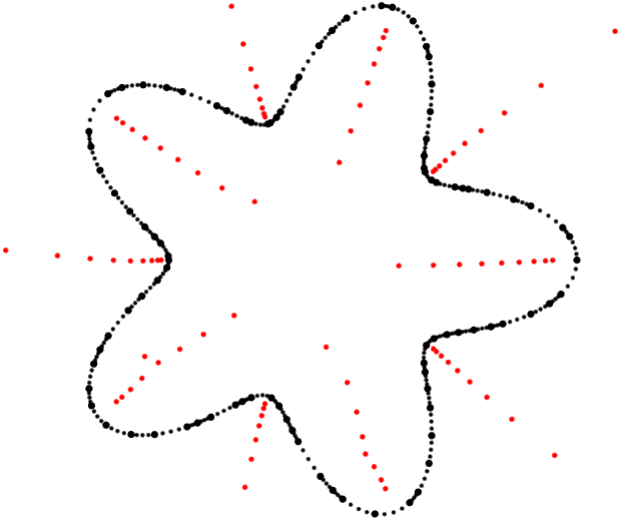}
		}
		\caption{Approximation of $v(z)=(\Re(z))^2$ on the parametric starfish curve $\Gamma(t) = (1+0.3\cos(5\pi\,t))e^{i\pi\,t}$ using the \textsc{Matlab} routine \texttt{aaazp}. Smaller black dots represent samples (381) generated on the fly, larger black dots are support points (95), and red dots are poles (94) of the rational approximation. The tolerance was loosened to $10^{-6}$ to clearly visualize the curve samples; accordingly, the computing time is reduced to half a second.}
		\label{fig:starfish-aaazp}
	\end{figure}

	To circumvent the need for manual grid refinement, a ``continuum'' variant of the AAA algorithm was introduced in \cite{Driscoll24}. Though less explored in applications (see for example \cite{Costa2024}), this version adopts an adaptive strategy, dynamically populating the sample set in regions where the error is greatest, thereby automatically ensuring the necessary clustering of both sample and support points according to more or less rapid variations in data and geometry. While existing implementations focus on the unit interval (\texttt{aaax}) and the unit circle (\texttt{aaaz}), the \texttt{aaazp} routine, an extension developed by the author, generalizes these adaptive capabilities to generic curves in the complex plane.
	For example, the \textsc{Matlab} lines
	\begin{verbatim}
		g = @(t) (1 + 0.3*cos(5*pi*t)).*exp(1i*pi*t);
		v = @(z) real(z).^2;
		r = aaazp(u,g,300,0,1e-13,1,0,2);
	\end{verbatim}
	approximate the data $v(z) = (\Re(z))^2$ on a starfish-shaped boundary $\Gamma(t)$ parametrized by the function handle \texttt{g}, returning a rational function \texttt{r} in barycentric form
	\footnote{
		Other parameters from third to last: 300 = maximum approximant degree; 0 = flag for activating Lawson iterations (off); 1e-13 = required tolerance; 1 = flag for indicating whether the function is meromorphic (on); 0 = flag for plotting (off); 2 = maximal degree of function derivative.
	}.
	This achieves a degree of 191 and an error of $7.7\cdot10^{-14}$ in approximately 3 seconds on a standard laptop, while also returning the first and second derivatives (Figure~\ref{fig:starfish-aaazp}). The function call strictly follows the syntax of \texttt{aaaz}, with differences only in the second, penultimate, and last parameters, representing the curve parametrization, the maximal degree of returned derivatives, and the initial guess for support points, respectively. Under the hood, code optimizations avoid redundant function evaluations at each step, greatly improving speed with respect to \texttt{aaaz}, and final accuracy is verified on a grid seven times finer than the original. Consequently, sample points automatically cluster as required in regions of strong curvature or rapidly varying boundary conditions. A key observation, the importance of which will become clearer in the following sections, is that the support points cluster in a similar fashion, and for a rational approximation of degree $J$ their number is exactly $J+1$.
		
\begin{figure}[!h]
		\centering{
			\includegraphics[width=0.9\linewidth]{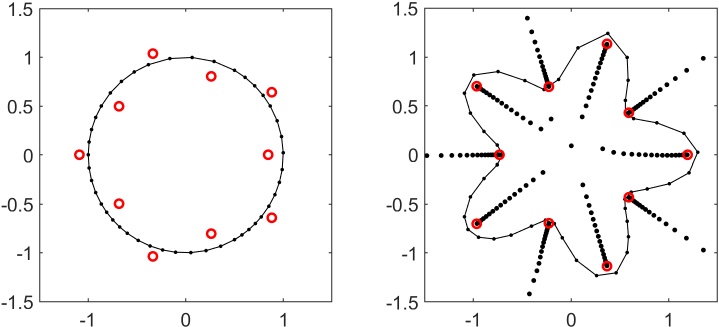}
		}
		\caption{Parametric mapping from the unit circle $C$ to the starfish curve $\Gamma$ of Figure (\ref{fig:starfish-aaazp}). On the left, red circles represent the zeros of $M'$. On the right, the starfish is shown with the branch cuts of the rational approximation of its Schwarz function (delineated by streams of black dots) and the images of the zeros of $M'$ through $M$ (red circles). With this particular parametrization, the zeros identify the branch points of $S(z)$. Only 45 sample points are necessary to reach a maximum error of $1.8\cdot10^{-15}$, explaining the seemingly coarse resolution of the curve.}
		\label{fig:disk2starfish}
	\end{figure}

	The same routine provides a convenient and efficient method for computing the Schwarz function
	\begin{verbatim}
		[Sf,polSf] = aaazp(@(z) conj(z),g,300,0,1e-13,1,0);
	\end{verbatim}
	along with the mapping $M$ from the unit circle $C$ in the $\zeta$-plane onto $\Gamma$ and the zeros \texttt{zerdM} of its derivative. This is facilitated by the algorithm's ability to differentiate a barycentric form, a feature adapted from standard AAA (\cite{Schneider86}), as shown in Figure~\ref{fig:disk2starfish}.
	\begin{verbatim}
		[M,~,~,~,~,~,Z] = aaazp(@(z) g(angle(z)/pi),@(t) exp(1i*pi*t),...
		  300,0,1e-13,1,0,1);
		[~,~,~,zerdM] = aaa(M{2}(Z),Z,'tol',1e-3);
	\end{verbatim}
	In this context, \texttt{Z} represents samples generated automatically on $C$. Thus, the continuum AAA algorithm for generic curves not only clusters an appropriate number of support points based on the characteristics of the data and the geometry, but also provides an efficient means of locating the singularities that limit analytic continuation.

	\section{The AAALS Algorithm for Helmholtz Scattering}\label{s:4}
	In this section we describe the procedure for solving exterior sound-soft Helmholtz scattering problems using the AAALS algorithm, providing extensive numerical results and comparisons. Our approach exploits to a great extent the fundamental results on the analytic continuation of the solution presented in \cite{BarnettBetcke08}, and in particular Conjecture 9 regarding the rate of convergence of the MFS for Helmholtz problems on general analytic domains.
		
	\subsection*{Theoretical foundation of the algorithm}
	The theoretical foundation for MFS point placement strategy is encapsulated in the following conjecture from \cite{BarnettBetcke08}, adapted using the notations of (\ref{eq:helmholtz}) and (\ref{eq:helmholtz-approx}) and section \ref{s:2}.
	\begin{uconj}
		Let $\varepsilon$ be the maximum error of the MFS measured on $\Gamma = \partial\Omega$, by placing the MFS points equally distributed in conformal angle at a conformal distance $R$ inside $\Gamma$. Let $\rho < 1$ be the conformal radius of the closest (in the sense of conformal radius) singularity of the analytic continuation of $u$. Then
		\begin{equation}
			\varepsilon \leq
			\begin{cases} 
				K \rho^{-J/2}, & \rho > R^2, \\ 
				K R^{-J}, & \rho < R^2, 
			\end{cases}
		\end{equation}
		where $K$ is a constant that may depend on $\Gamma$, $k$, $R$ and $v$, but not $J$.
	\end{uconj}
	To leverage these results effectively, we perform the rational approximation on the unit circle $C$ in the computational $\zeta$-plane rather than on the physical boundary $\Gamma$ in the $z$-plane. This choice is driven by three key factors:
	\begin{enumerate}
		\item \emph{Conformal geometry and MFS placement:} The convergence of the MFS is governed by the relative positions of singularities and sources in terms of \emph{conformal radius}. Working in the computational $\zeta$-plane allows us to precisely determine the conformal radius $\rho$ of the singularities detected by AAA and, consequently, to place the MFS source points at a specific conformal distance $R$. Achieving this specific geometric arrangement directly on $\Gamma$ without the explicit mapping would be practically impossible.
		\item \emph{Precision and robustness:} Analytic continuation is invariant under conformal transplantation, much like the singularities themselves. Approximating directly on $\Gamma$ and then mapping the resulting poles to the $\zeta$-plane via $m(z) = M^{-1}(z)$ would introduce significant numerical challenges, specifically the generation of ``spurious'' approximation poles that require complex filtering. Moreover, identifying which singularities lie outside the conformality region remains non-trivial. Working directly on $C$ circumvents these issues entirely.
		\item \emph{Computational efficiency:} The use of a standardized domain allows for the efficient recycling of support points between sequential calls to the solver, leading to a substantial reduction in the overall computational overhead.
	\end{enumerate}
		
	\subsection*{Practical aspects}
	To achieve this, we deviate from the standard practice of working directly on the physical boundary $\Gamma$, which is typical of AAALS applications for Laplace problems. Instead, we first perform a rational approximation of the real part of the complex function defining the incident wave $u_{inc}(z)$ on the unit circle $C$. The approximation is carried out using the continuum AAA variant (\texttt{aaazp}). The use of the continuum mode is not merely a convenience but a necessity: static discretizations often fail to capture the rapid variations of the boundary data caused by singularities located near the boundary. The adaptive nature of \texttt{aaazp} ensures that the sampling density automatically matches the local complexity of the function and the geometry.
		
	Finally, we address the placement of the MFS sources. Standard AAALS typically uses the poles of the rational approximant directly as sources. However, following the analysis in \cite{BarnettBetcke08}, we adopt a more robust strategy. We use the poles identified by AAA, augmented with the zeros of the mapping derivative $M'$, to define a ``shielding curve'' $\gamma$ within the conformality region $A$ interior to $C$. We then select points $\{\pi_j\}$ at well-defined distances between $\gamma$ and $C$, whose radial projections onto the latter correspond to the support points of the approximation. When mapped to the physical domain, the images $\{p_j = M(\pi_j)\}$ naturally cluster and spread according to the variations in the boundary data and the gradient of the mapping function. This results in a singularity-adapted distribution of sources that is specifically tailored for the MFS, providing a robust approximation basis superior to raw pole reuse for low and high wave numbers alike.
		
	\subsection*{Implementation}
	Below we outline the steps of the algorithm in greater detail. All approximations are computed using the \texttt{aaazp} routine to exploit its code optimizations.
	\begin{enumerate}
		\setlength{\itemsep}{5pt}
		\setlength{\parskip}{0pt}
		\setlength{\parsep}{0pt}
		\item Define a parametrization for $\Gamma(t), \; t \in [-1,1]$ in the $z$-plane, and a complex boundary condition $u_{inc}(z)$ expressing the values to be canceled by the outward-radiating sources (the sound-soft condition).
		\item Approximate the mapping $M: C = e^{i\pi t} \rightarrow \Gamma(t)$ and compute the set $\mathcal{O}$ of zeros of $M'$ interior to $C$ in the $\zeta$-plane. \label{it:appr-M}
		\item Approximate $v = \Re(u_{inc}(\Gamma(t)))$ with $t = \arg(\vartheta)/\pi$ on $C = e^{i\vartheta}, \; \vartheta \in [-\pi,\pi]$ to identify the set $\mathcal{Q}$ of poles interior to $C$ that limit its analytic extension. Notably, this approach remains effective even for non star-shaped domains. \label{it:appr-C}
		\item Determine a conforming polygon $P$ from the set consisting of $\mathcal{Q \cup O}$, utilizing, for instance, the \textsc{Matlab} \texttt{boundary} routine with a low or moderate shrink factor (typically between 0.1 and 0.5).
		\item Interpolate a curve $\gamma(\vartheta)$ from the vertices of $P$; the use of \texttt{interp1} (specifically with \texttt{pchip} or \texttt{makima} methods) is recommended and consistently outperforms \texttt{spline} interpolation, which tends to introduce unwanted overshoots and ringing artifacts.
		\item Approximate $v$ again, but this time directly on $\Gamma$, to determine $J$ support points $\{s_j\}$ and $N \gg J$ boundary samples $\{z_n\}$. The support points $\{s_j\}$ on $\Gamma$ map to $\{e^{i\theta_j}\}$ on $C$. \label{it:appr-G}
		\item Compute the points $\{\gamma(\vartheta_j) = \rho_j\,e^{i\theta_j}\}$ on $\gamma$, and subsequently the points $\{\pi_j = \sqrt{\rho_j}\,e^{i\theta_j}\}$ located at an optimal distance from $C$.
		\item Compute the final positions of the scattering sources $\{p_j = M(\pi_j)\}$ inside $\Gamma$.
	\end{enumerate}
	A refined detail in these approximation steps is that we can initialize each stage using the support points from the previous approximation. Specifically, the support points resulting from step~\ref{it:appr-M} are fed into step~\ref{it:appr-C}, and those from step~\ref{it:appr-C} are used in step~\ref{it:appr-G}. This strategy leads to a dramatic reduction in total computing time. Furthermore, experimental evidence shows that very few additional support points are appended during step~\ref{it:appr-G}, especially as the wavenumber increases.
		
	Ultimately, the scattered wave $u_s$ in the exterior of $\Gamma$ is computed via a least squares fit to the incident wave $u_{inc}$ using the sources $\{p_j\}$, while the total solution outside $\Gamma$ is found by subtracting $u_{inc}$ from $u_s$. We utilize the following multipole expansion:
	\begin{equation}
		u_s(z) \approx \sum_{j=1}^{J} \sum_{r=-R}^R c_{j,r} H_{|r|}^{(1)}(k|z-p_j|) \left( \frac{z-p_j}{|z-p_j|} \right)^r,
		\label{eqn:expansion}
	\end{equation}
	where $H_r^{(1)}$ are the Hankel functions of the first kind (outward-radiating) and $\{c_{j,r}\}$ are the coefficients resulting from the least squares problem. For the sake of conciseness, let
	\[
	\mathbf{H}(r,p_j,z_n) = \begin{pmatrix} H_{|r|}^{(1)}(k|z_n-p_j|)\left(\frac{z_n-p_j}{|z_n-p_j|}\right)^r \end{pmatrix} \in \mathbb{C}^{N\times(2JR+J)}
	\]
	be the matrix constructed with the parameters varying as described, where $[z_1 \dots z_N]^T$ is a column vector and $[p_1 \dots p_J]$ is a row vector. This leads to the following system:
	\begin{equation}
		\mathbf{H}\cdot
		\begin{pmatrix}
			[c_{1,-R}\; ...\; c_{N,-R}]^T\\
			[c_{1,-R+1}\; ...\; c_{N,-R+1}]^T\\
			\vdots \\
			[c_{1,R}\; ...\; c_{N,R}]^T
		\end{pmatrix} = 
		\begin{pmatrix}
			u_{inc}(z_1)\\
			\vdots\\
			u_{inc}(z_N)
		\end{pmatrix}
		\label{eqn:LS}
	\end{equation}
	which constitutes an overdetermined $N \times (2JR+J)$ system in $2JR+J$ unknowns, solved using the \textsc{Matlab} backslash operator.
	\begin{figure}[!h]
		\centering{
			\includegraphics[width=0.8\linewidth]{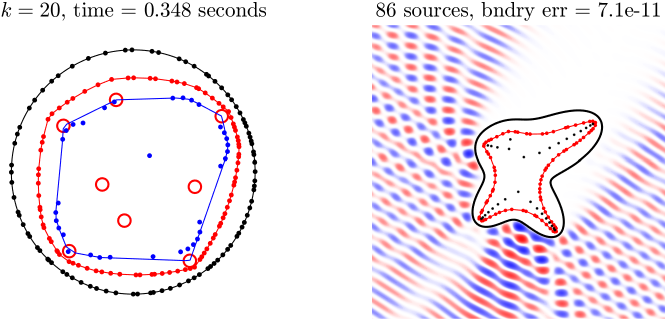}\\
			\medskip
			\includegraphics[width=0.8\linewidth]{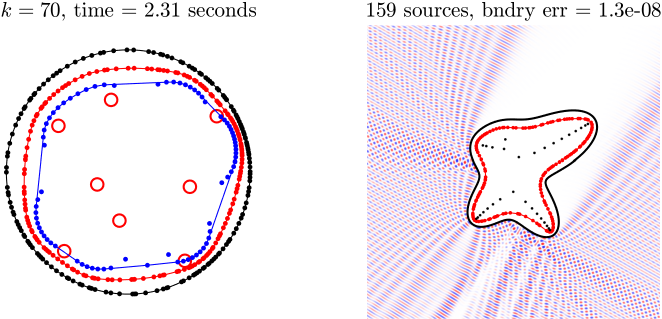}\\
			\medskip
			\includegraphics[width=0.8\linewidth]{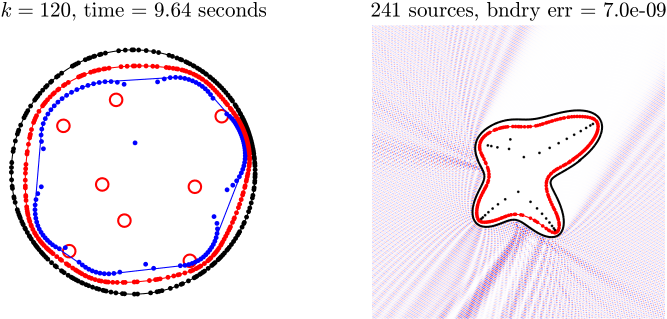}
		}
		\caption{``Sound-soft'' scattering problems for incident plane waves with angle $\pi/3$ and different wave numbers. On the left, the blue line $\gamma$ bounds the zeros of $M'$ (red circles) and the poles (blue dots) of the rational approximation of $v = \Re(u_{inc}(\Gamma))$ interior to $C$; black dots are support points of the rational approximation and red dots are their radial projections at an optimal distance between $\gamma$ and $C$. On the right, the solutions; red dots are scattering sources and black dots delineate the branch structure of the Schwarz function.}
		\label{fig:rand_wave}
	\end{figure}
	\begin{figure}[!t]
		\centering{
			\includegraphics[width=0.8\linewidth]{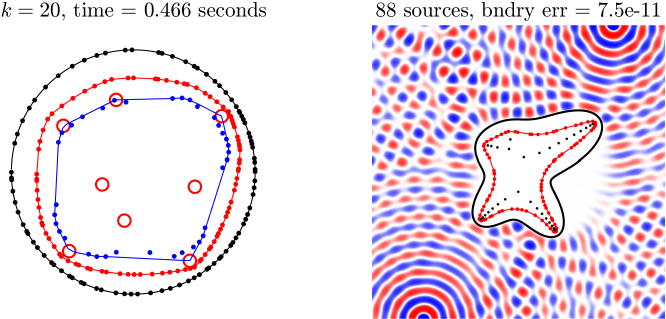}\\
			\bigskip
			\includegraphics[width=0.8\linewidth]{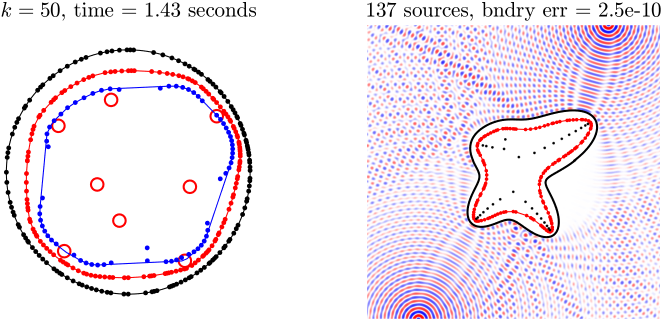}\\
			\bigskip
			\includegraphics[width=0.8\linewidth]{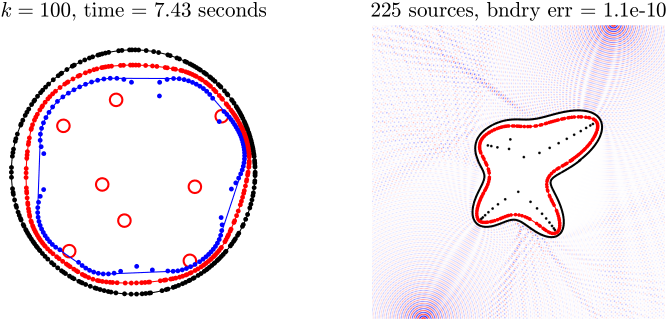}\\
		}
		\caption{The same curve $\Gamma$ of figure \ref{fig:rand_wave} for two point sources of equal wave number.}
		\label{fig:rand_src}
	\end{figure}
	
	\subsection*{Numerical results}By employing the parameter values specified as follows, the algorithm consistently achieves 6 or more digits of accuracy (typically more), as measured by the $L^\infty$-norm of the boundary condition on a grid four times finer than the original, taking times of the order of seconds on a laptop. Although some of these parameter settings represent engineering choices, they have proven to be robust across a large number of randomly generated smooth Jordan curves for wavenumbers up to $k=150$ and beyond. The \texttt{aaazp} tolerance for approximating $v$ on the circle is set to $10^{-6}$, while a tolerance of $10^{-8}$ is used for approximating $M$. While we limit the expansion in (\ref{eqn:expansion}) to $R = 2$, it is worth noting that including higher-order Hankel functions requires a higher sampling density (9 samples instead of the standard 3) between support points for an accurate least squares solution. Computationally, including Hankel functions of order 2 and increasing the number of sample points is far more efficient than simply tightening the approximation tolerance. Additionally, the matrix columns in (\ref{eqn:LS}) are scaled by the distance of each $\{p_j\}$ from the boundary, an essential step for preserving the condition number as sources approach the boundary. This ensures that the Hankel functions attain a maximum value of $1$ on the curve (this key adjustment was suggested by Andrea Moiola).
		
	Figures~\ref{fig:rand_wave} and \ref{fig:rand_src} display the results for increasing values of $k$, considering an incident plane wave traveling at $\pi/3$ and two distinct point sources, respectively. In the left columns, it is evident that for low wavenumbers, the curve $\gamma$ is primarily determined by the zeros of $M'$, while the poles of the boundary condition approximation fall closer to the center, thus outside the domain of conformity. As $k$ increases, this situation is reversed.
		
	\subsection*{Validation} 
\begin{figure}[!t]
		\centering{
			\includegraphics[width=0.42\linewidth]{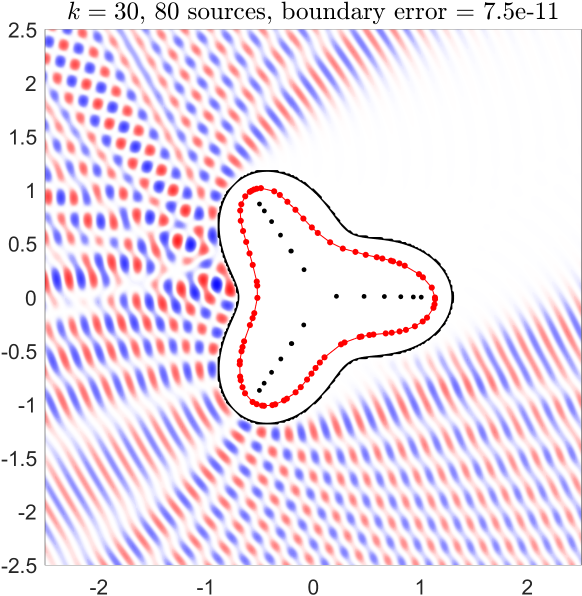}\qquad
			\includegraphics[width=0.5\linewidth]{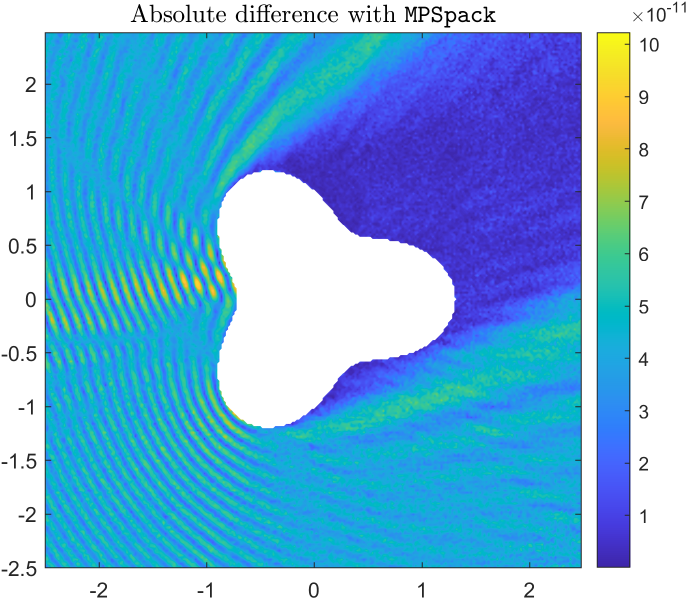}
		}
		\caption{On the left, the real part of the solution computed with the proposed algorithm for the sound-soft scattering example found in \cite{MPSpackTutorial}. On the right, the absolute point-wise difference with respect to values provided by \texttt{MPSpack}; deviation takes into account both real and imaginary parts of the waves, and is about $10^{-10}$ at its peak.}
		\label{fig:comparison}
	\end{figure}

The algorithm is relatively straightforward to implement, with a significant portion of the computational effort handled autonomously by the continuum variant of the AAA algorithm, which provides mapping functions, singularities, support points, and boundary samples in just a few calls. One might question the accuracy of these results, which in Figures~\ref{fig:rand_wave} and \ref{fig:rand_src} is measured specifically on the boundary: are these digits representative of the solution's quality elsewhere? To address this, we use the \texttt{MPSpack} toolbox for \textsc{Matlab} (\cite{MPSpackManual}) as a reference, which is specifically designed for high-accuracy solutions using basis functions. We solve the sound-soft scattering problem provided as an example in \cite{MPSpackTutorial}, featuring a plane wave with $k = 30$ incident at an angle of $\pi/6$ on a trefoil domain. While \texttt{MPSpack} requires a basis of 210 functions and 250 quadrature points to achieve an $L^\infty$ error of $8 \cdot 10^{-12}$ on $\Gamma$, our method automatically generates 80 sources and 720 boundary samples, achieving a comparable error of $7.5 \cdot 10^{-11}$ in a fraction of a second. The resulting solution and the point-wise absolute difference compared to \texttt{MPSpack} are shown in Figure~\ref{fig:comparison}.

	\section{Concave, Elongated, and Non-Analytic Curves}\label{s:5}
	In this section, we summarize several results for more challenging geometries, specifically those involving concavities or strong elongations that introduce the numerical ``crowding problem'' of conformal mapping into the analysis. Furthermore, while the AAALS-Helmholtz algorithm is primarily designed for analytic curves, it demonstrates a robust capability to handle corner singularities.
	\begin{figure}[!t]
		\centering{
			\includegraphics[width=0.45\linewidth]{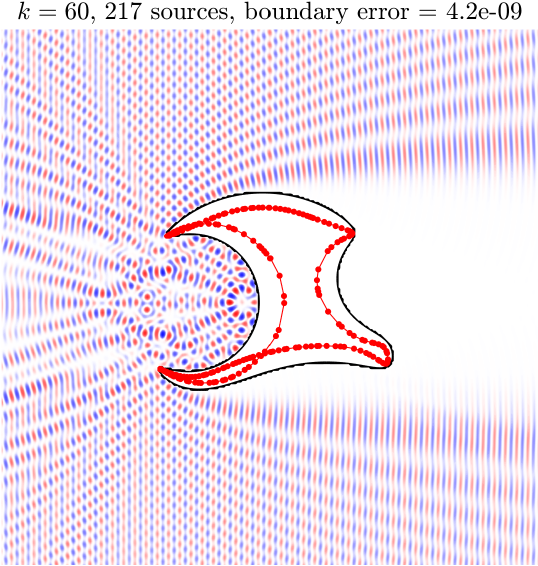}\qquad
			\includegraphics[width=0.45\linewidth]{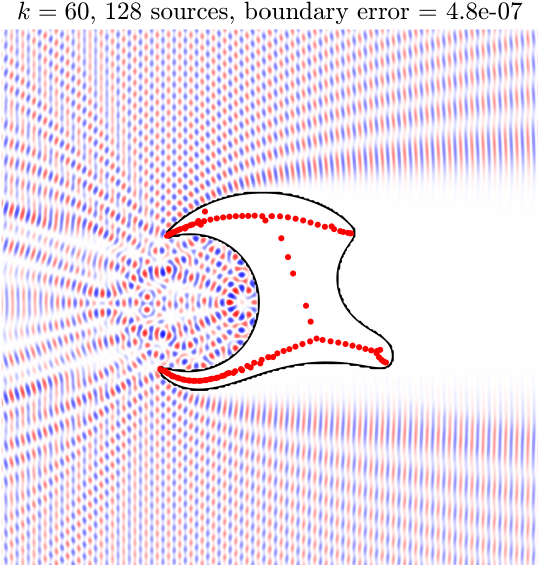}\\
			\bigskip
			\includegraphics[width=0.45\linewidth]{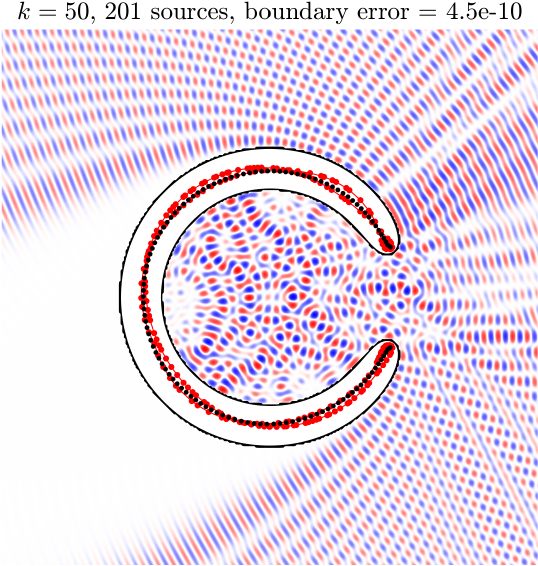}\qquad
			\includegraphics[width=0.45\linewidth]{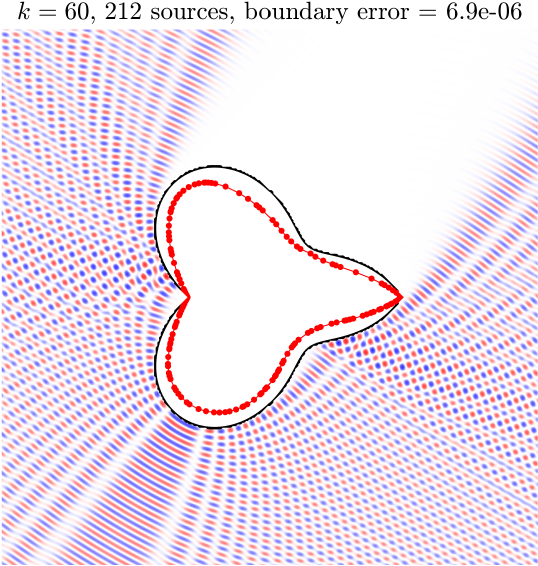}\\
		}
		\caption{Four examples of scattering problems involving numerically challenging curves. Clockwise from top-left: the generalized crescent from [] solved with the proposed algorithm; the same crescent solved for comparison using source points determined by a straightforward AAA approximation of $v = \Re(u_{inc}(\Gamma))$; a corral with sources (red) strictly following the branch structure of the Schwarz function (black); and a ``sharktooth'' featuring two corner singularities. While the application to curves with corner singularities remains largely unexplored, these results are encouraging.}
		\label{fig:crescents}
	\end{figure}
		
	Various examples are illustrated in Figure~\ref{fig:crescents}. In the upper-left panel, the algorithm's performance is demonstrated on the generalized crescent curve from \cite{BarnettBetcke08}:
	\[
	Z_{GC}(z) = e^{iz}-\frac{0.1}{e^{iz}+0.9}-\frac{0.07+0.02i}{e^{iz}-0.8-0.2i}+\frac{0.2}{e^{iz}-0.2+0.5i}.
	\]
	Although accurate results can be achieved for significantly higher wavenumbers, the choice of $k=60$ facilitates a direct comparison with the standard computation using poles obtained from a straightforward AAA approximation of $\Re(u_{inc}(\Gamma))$, shown in the upper-right panel. It appears that for geometries characterized by narrow elongations, the branch structure identified by the AAA algorithm is sufficient for high-quality approximations. Indeed, an inspection of the unit circle in the $\zeta$-plane clarifies that the curve $\gamma$ is predominantly determined by the zeros of $M'$, even for $k \gg 20$. This behavior is further confirmed by the ``corral'' shape in the lower-left panel, where the point sources (in red) strictly follow the branch structure of the Schwarz function (in black). This structure, in turn, matches that obtained by the AAA approximation of the real boundary data; consequently, as expected, the simpler computational method would yield nearly identical results in this case.
		
	Finally, a ``sharktooth'' curve proposed by Nick Trefethen is presented in the lower-right panel. The algorithm performs reasonably well, although the maximum achieved accuracy was limited to approximately $10^{-6}$. We believe this ultimately stems from the reduced performance of the \texttt{aaazp} routine in the presence of corner singularities, specifically regarding adequate sampling, though this issue has not yet been investigated in depth.

	\section{Application to Laplace Problems}\label{s:6}
	The first application of AAALS to the solution of PDEs was introduced for Laplace problems in \cite{Costa20} and \cite{CostaTrefethen23}:
	\[\Delta f = 0 \quad \text{in } \Omega, \quad f(z) = v(z) \quad \text{on } \Gamma = \partial\Omega. \]
	Even so, in over five years, not much formal theory has been established to explain why it works so effectively, especially regarding the ``double poles'' technique for improving (effectively doubling) the accuracy compared to standard single-pole computations. A discussion of this can be found in section 24 of \cite{NakatsukasaTrefethen26}, and the justification, based on standard arguments of potential theory, may be outlined as follows. Essentially, geometry-driven problems possess the same, or nearly the same, branch structure as the Schwarz function, which is delineated by streams of poles. Those returned by the AAA algorithm are well-placed to approximate a real-valued function $v(z)$ on $\Gamma$; we exploit its branch structure along with the fact that if a degree $n$ rational approximation $r(z)$ interpolates a function $w(z) = v(z) + ih(z)$, with $h$ the harmonic conjugate of $v$, at $2n+1$ points $\{s_k\}$, then the Hermite integral:
	\[w(z) - r(z) = \frac{1}{2\pi i} \int_\Gamma \frac{\varphi(z)}{\varphi(t)} \frac{w(t)}{t-z} \, dt \]
	features a squared denominator:
	\[\varphi(z) = \prod_{k=0}^{2n}(z-s_k) \Bigg/ \prod_{k=1}^{n}(z-p_k)^2. \]
	Therefore, we argue that after discarding a subset of poles from the AAA rational approximation for $v$, by doubling the remaining $\{p_k\}$, we obtain, at least asymptotically, the same accuracy provided by a nearly-optimal set for $f$ (and for $w$).
		
\begin{figure}[!t]
		\centering{
			\includegraphics[width=0.8\linewidth]{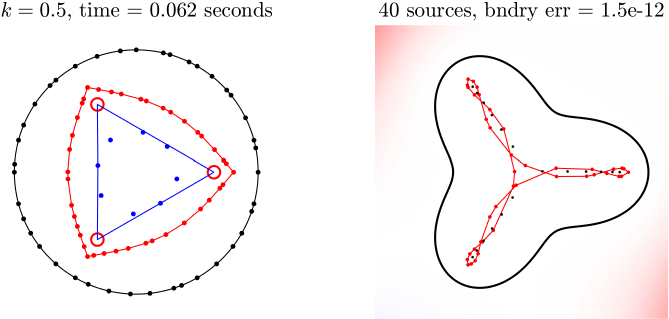}\\
		}
		\caption{Result of the algorithm of section \ref{s:4} on a trefoil $\Gamma$ for an incident plane wave of number close to zero. Red dots between $\gamma$ (blue triangle) and $C$ on the left are mapped through $M$ onto source points on the right, verging on interior branch cuts of the Schwarz function for $\Gamma$.}
		\label{fig:lowk}
	\end{figure}

	We propose an alternative justification based on results for the MFS. First, since Laplace problems can be viewed as the limit of Helmholtz problems as $k \to 0$, we can examine what occurs for low wavenumbers using the algorithm described in section~\ref{s:4}. As shown in Figure~\ref{fig:lowk}, for $k = 0.5$, the source points follow closely the branch cuts of the Schwarz function for $\Gamma$, which in turn behave like those of the approximant to the real part of the incident wave as the problem becomes purely geometry-driven. This behavior allows for a straightforward understanding of why it has been observed that accurate solutions to Helmholtz problems can be obtained directly from the poles returned by AAA approximation only for low wavenumbers, i.e., up to approximately $k=20$ (\cite{Gopal19,Gopal21,Ginn23}).
		
	The key result allowing for a deeper understanding comes from a paper by Katsurada and Okamoto (\cite{Katsurada96}) on the collocation points of the MFS for potential problems, where they specifically address interior Laplace problems. Their method requires selecting $N$ collocation points and $N$ charge points, representing the singularities of the fundamental solution. The former are freely selected from the curve $\Gamma$, while charge points are subsequently determined through a conformal map $\Psi$, defined in an annular neighborhood of $C$ and mapping $C$ onto $\Gamma$, computed via complex FFT. Eventually, a least squares problem is solved to determine the coefficients for approximating the solution $f$ in $\Omega$. This method yields effective charge points, providing an exponential decrease in the boundary error, only when provided with good collocation points. While \cite{Katsurada90} and \cite{Katsurada94} prove their theoretical existence, the authors note the lack of a precise rule for choosing them or an appropriate value for the radius $R > 1$ for placing charge points around $C$. Moreover, accurate results are contingent upon the precise computation of $\Psi$ via FFT.
		
	If we reconsider the algorithm of section~\ref{s:4} in light of these results, we have the mapping $M$ in place of $\Psi$, while the curve $\gamma$ (or $\delta$) improves the choice of a fixed distance from $C$ to place charges for exterior (or interior) problems. Meanwhile, the $J+1$ support points of the rational approximant of degree $J$ to the real-valued function $v$ on $\Gamma$ serve as collocation points where the error is zero. Remarkably, the AAA algorithm alone is capable of providing all this critical information.
	\begin{figure}[!t]
		\centering{
			\includegraphics[width=0.48\linewidth]{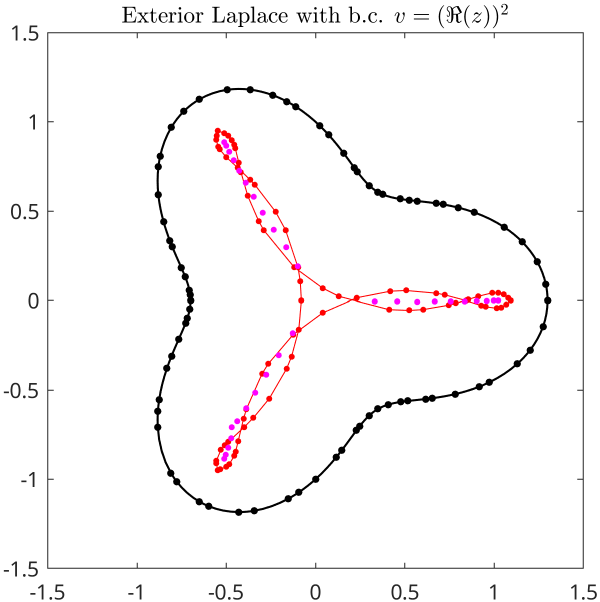}\quad
			\includegraphics[width=0.48\linewidth]{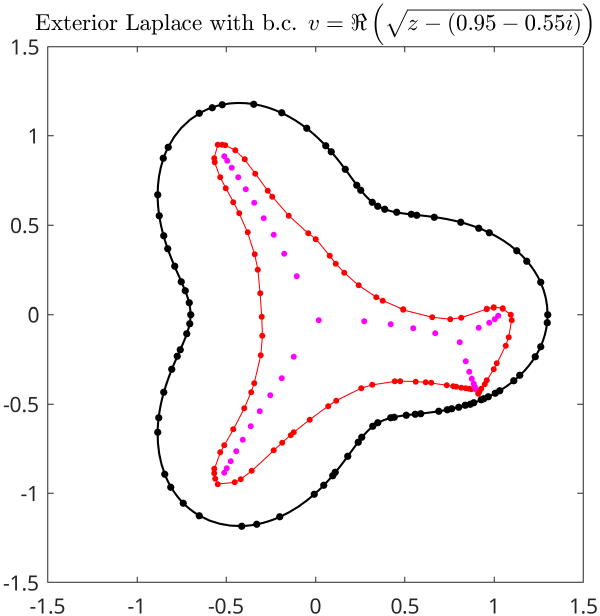}
		}
		\caption{Two Laplace problems for the exterior to $\Gamma$. On the left, the problem is completely geometry-driven; poles (red dots) are derived from support points (in black on the curve) and are the same in number. They are positioned by our algorithm about branch cuts delineated by those returned by the AAA approximation of the boundary condition (in magenta), and notably their alignment oscillates between different Riemann sheets. The original AAA accuracy of $10^{-12}$ is preserved. On the right, the boundary condition exhibits a singularity close to $\Gamma$; accuracy is not fully restored due to the difficulty to interpolate a close-fit line $\gamma$, hence poles don't follow branch cuts closely as they should.}
		\label{fig:lapl_ext}
	\end{figure}
	Two examples illustrate this. In Figure~\ref{fig:lapl_ext} (left), we solve for $v(z) = (\Re(z))^2$ in the exterior. The \texttt{aaazp} routine is used on $\Gamma$ with a tolerance of $10^{-12}$ to efficiently find 32 interior poles (magenta dots) and 74 support points $\{s_j\}$ (black dots on the curve). The images of these support points on the unit circle in the $\zeta$-plane are adequately repositioned between $C$ and $\gamma$, and eventually mapped back via $M$ to find the poles $\{p_k\}$ for the approximant to $f$, which match the number of support points. These poles end up placed closely around the magenta dots; their trajectory also jumps between different sides of the branch cuts, which is expected when approaching them closely. The final boundary error is approximately $10^{-13}$. For geometry-driven problems, branch points are determined by the Schwarz function (see Figure~\ref{fig:lowk}, left), and small deviations in the shape of branch cuts have little influence on final accuracy. More significant is the fact that poles (the images of support points) tend to merge there, and since they are exactly $J+1$, they preserve accuracy. From a practical perspective, simply doubling poles in the interior is often the most efficient method, even if they are unevenly distributed between interior and exterior. In Figure~\ref{fig:lapl_ext} (right), the boundary condition is changed to $v(z) = \Re\left(\sqrt{z-(0.95-0.55i)}\right)$, introducing a singularity close to $\Gamma$. Our method yields a final error of about $10^{-10}$; this imperfect result is due to the fact that $\gamma$ does not bound the domain of conformity efficiently and to its near maximal extension, which keeps poles at a distance from the branch cuts. Again, pole doubling proves more efficient and consistent, resulting in an error of approximately $10^{-14}$.
		
\begin{figure}[!t]
		\centering{
			\includegraphics[width=0.495\linewidth]{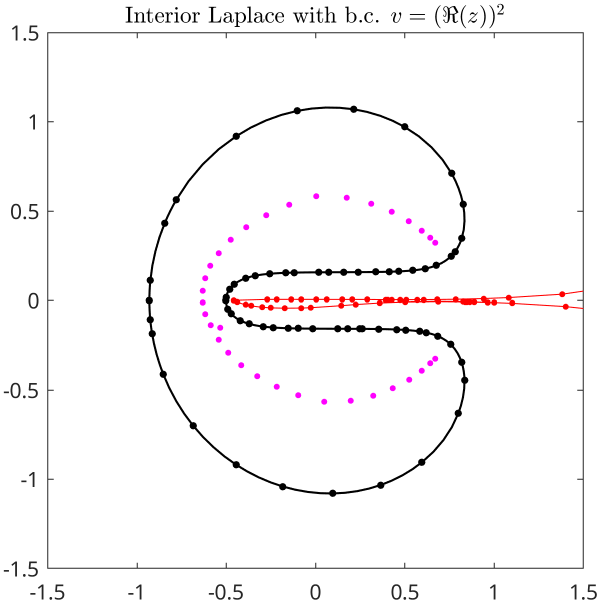}\quad
			\includegraphics[width=0.47\linewidth]{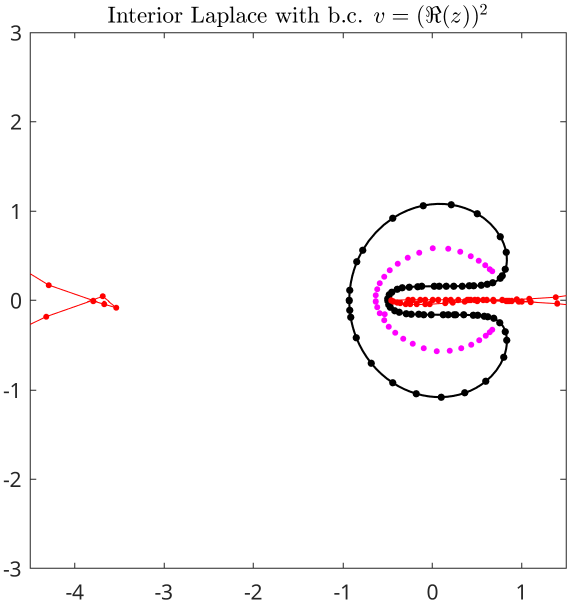}
		}
		\caption{An interior Laplace problem for a ``Pac-man'' curve $\Gamma$. On the left, the line with poles (in red) derived by the algorithm of section\ref{s:4} from support points, surrounding those returned by AAA approximation (magenta); the latter are absent near the leftmost curved side. On the right, a zoom out showing also poles derived from the leftmost support points on $\Gamma$; other poles surround the curve farther and are not shown for clarity.}
		\label{fig:pacman}
	\end{figure}

	Finally, we can explain the failure of AAALS in the presence of certain symmetries. Consider the ``Pac-man'' geometry proposed by Nick Trefethen in Figure~\ref{fig:pacman}, with a slightly shortened inlet to make the problem numerically tractable. The AAA algorithm for both $v(z) = (\Re(z))^2$ and the Schwarz function fails to place exterior poles near the circular side on the left of the curve, which has severe consequences for the accuracy of the approximant to $f$ in the interior, even with pole doubling. However, we know they should be present, as artificially adding them and performing a subsequent ``AAA compression'' confirms their existence. The proposed algorithm, by accounting for support points that appear along the entire curve, places several poles in the correct locations, and the final accuracy increases from 3 digits (standard AAALS) to 5. While the original tolerance of $10^{-6}$ is not fully restored due to the extreme geometry, the principle remains valid. Other similar cases in the literature, such as the channel with a deep inlet in \cite{Xue24}, can be explained accordingly. We have not yet identified a computational shortcut to make such ``ghost'' poles materialize more easily.

	\section{Discussion}\label{s:7}
	In this paper, we have presented a comprehensive numerical framework for the solution of exterior sound-soft Helmholtz scattering problems with high wave numbers by integrating the Method of Fundamental Solutions (MFS) with the AAA rational approximation algorithm. The resulting AAALS-Helmholtz method explicitly addresses the lack of formalized literature regarding adaptive rational approximation in wave propagation. Framing the MFS within the context of rational approximation, we have validated the hypothesis that an automated strategy for singularity placement, dictated by the analytic continuation of boundary data, can vastly outperform traditional heuristics.
		
	A cornerstone of this approach is the introduction of continuum AAA. As conjectured in previous studies on analytic domains, the success of the MFS relies on a choice of source curve that avoids the singularities of the boundary data's analytic continuation. The \texttt{aaazp} algorithm automates this task, creating singularity-adapted distributions that ensure exponential convergence. Furthermore, we have provided a concrete answer to the open question regarding accuracy enhancement: our experiments demonstrate that augmenting the basis with low-degree Hankel series (specifically up to order $R=2$), coupled with density-scaled sampling, is a highly effective strategy to reach high precision without resorting to ill-conditioned bases or artificially increased forcing frequencies.
		
	Theoretical connections established in this work also clarify the double poles technique for Laplace problems. We have shown that the source placement strategy derived for Helmholtz naturally reduces to the optimal dipole placement for potential problems as $k \to 0$, providing a justification for what was previously an empirical observation.
		
	While the solution of the linear system is rapid due to the small basis size achieved by the AAALS approach, the evaluation of the solution in the bulk for high-resolution visualization remains computationally intensive. A natural next step for large-scale applications would be the integration of this solver with a Fast Multipole Method (FMM) to accelerate the evaluation phase. Other research can also target Neumann and transmission problems, where adaptive singularity tracking could prove even more decisive. In conclusion, the AAALS-Helmholtz method is a promising alternative among meshless solvers, offering a robust, adaptive, and theoretically grounded tool for engineering computations.
		
	\subsection*{Acknowledgements}
	I am grateful to Andrea Moiola and Nick Trefethen for several helpful discussions and for their comments on earlier versions of this work.
	
	\bigskip
	\bibliographystyle{amsalpha}
	\bibliography{x}
	
\end{document}